\documentclass{amsproc}

\numberwithin{equation}{section}
\newcommand{\C}{\mathbb{C}}
\newcommand{\R}{\mathbb{R}}
\newcommand{\Q}{\mathbb{Q}}
\newcommand{\Z}{\mathbb{Z}}
\newcommand{\To}{\rightarrow}
\newcommand{\cal}[1]{\mathcal{#1}}
\theoremstyle{plain}
\newtheorem{Thm}{Theorem}[section]

\newtheorem{Lem}[Thm]{Lemma}
\newtheorem{Prop}[Thm]{Proposition}

\theoremstyle{remark}
\newtheorem{Ex}{Example}
\newtheorem{remark}{Remark}
\newtheorem{definition}{Definition}
\newtheorem{Conj}{Conjecture}
\begin{document}
\title
{Unipotent Jacobian Matrices and Univalent Maps}
\author{L. Andrew Campbell}
\address{The Aerospace Corp. \\
M1-102 PO Box 92957 \\
Los Angeles CA 90009}
\email{campbell@aero.org}
\dedicatory{Combinatorial and Computational Algebra, Hong Kong, May 24--29, 1999}
\begin{abstract}
The Jacobian Conjecture would follow if it were known
that real polynomial maps with a unipotent Jacobian
matrix are injective. The conjecture that this is
true even for $C^1$ maps is explored here.
Some results known in the polynomial case are
extended to the $C^1$ context, and some special cases are
resolved.
\end{abstract}
\maketitle
\section{Introduction}
The focus of this paper is the unipotence
(all eigenvalues are 1)  
of the Jacobian matrix 
of a $C^1$ (continuously differentiable) map
from $\R^n$ to itself
and whether this implies its univalence (injectivity) or invertibility (bijectivity).
In the case of polynomial maps, unipotence is central to
reformulations of the Jacobian Conjecture. 
The paper is organized as follows:
a review of a number of important examples;
a comparison
of several established conjectures related to unipotence; 
a description of the goals of this paper;
a number of results that parallel what is known
in the polynomial case;
and finally, some partial results in the general $C^1$ context.
\section{Examples}
Consider some simple examples of maps with unipotent
Jacobian matrices and their explicit inverses.

\begin{Ex}
\label{Ex1}
Let $f: \R^2 \To \R^2$ be defined by
$$ (u,v) = f(x,y) = (x + 5\cos(3x+5y),\; y - 3\cos(3x+5y)) $$
Then the Jacobian matrix of $f$ is
$$
J(f) =
\begin{bmatrix}
\partial u/\partial x & \partial u/\partial y \\
\partial v/\partial x & \partial v/\partial y
\end{bmatrix}
=
\begin{bmatrix}
1 -15\sin(3x+5y) & -25\sin(3x+5y) \\
9\sin(3x+5y) & 1 + 15\sin(3x+5y)
\end{bmatrix}
$$
So $J(f)$ is unipotent; that is
$$
J(f)
=
\begin{bmatrix}
1 & 0 \\
0 & 1
\end{bmatrix}
+
\sin(3x+5y)
\begin{bmatrix}
-15 & -25 \\
9   &  15
\end{bmatrix}
$$
is the sum $I + N$, where $I$ is the ($2 \times 2$) identity
matrix and $N$ is a nilpotent matrix (a power of $N$ is $0$;
here it is easy to verify that $N^2 = 0$).
From the fact that $3u+5v = 3x+5y$, it follows easily
that the map $f$ is invertible, with inverse
$$ (x,y) = g(u,v) = f^{-1}(u,v) =
(u - 5\cos(3u+5v),\; v + 3\cos(3u+5v))
$$
\end{Ex}

\begin{Ex}
\label{Ex2}
Let $f: \R^4 \To \R^4$ be a $C^1$ upper triangular map:
$$ (s,t,u,v) = f(x,y,z,w) = (x + a(y,z,w), y + b(z,w), z + c(w), w + d) $$
where $a,b,c$ are $C^1$ functions of the indicated
variables and $d$ is a constant. 
$J(f)$ is unipotent; it is the sum of the ($4 \times 4)$ identity
matrix $I$ and a matrix $N$ that is strictly upper triangular
and hence nilpotent ($N^4 = 0$). The inverse of $f$ is
$$  (x,y,z,w) = f^{-1}(s,t,u,v) = $$
$$
(s - a(t-b(u-c(v-d),v-d),u-c(v-d),v-d),
 t -     b(u-c(v-d),v-d),
 u -         c(v-d),
 v -             d)
$$
\end{Ex}

\begin{Ex}
\label{Ex3}
Let $f: \R^3 \To \R^3$ be given by
$$ (u,v,w) = f(x,y,z) = (x + z\phi(x+zy), y - \phi(x+zy), z) $$
where $\phi$ is a $C^1$ function of a single variable. Then
$$
J(f) =
I +
\begin{bmatrix}
z\phi'(x+zy) & z^2 \phi'(x+zy) & \phi(x+zy) + zy\phi'(x+zy) \\
-\phi'(x+zy)   &  -z\phi'(x+zy)    &  -y\phi'(x+zy) \\
0          &     0        & 0
\end{bmatrix}
$$
which represents $J(f)$ as $I+N$, with $N$ nilpotent
($N$ is nilpotent since its upper left
$2 \times 2$ block is). 
From the fact that $u+vw = x+zy$, it follows, as in
the first example, that the map $f$ is invertible, with inverse
$$ (x,y,z) = g(u,v,w) = f^{-1}(u,v,w) =
(u - w\phi(u+vw), v + \phi(u+vw), w)
$$
\end{Ex}
\begin{Ex}
\label{Ex4}
Let $f: \R^3 \To \R^3$ be given by
$$ (u,v,w) = f(x,y,z) = (x+\phi(y-x^2),y+z+2x\phi(y-x^2),z-(\phi(y-x^2))^2)$$
where $\phi$ is a $C^1$ function of a single variable. Then
$$ J(f) = I +
\begin{bmatrix}
-2x\phi'(y-x^2) & \phi'(y-x^2) & 0 \\
2\phi(y-x^2) -4x^2\phi'(y-x^2) & 2x\phi'(y-x^2) & 1 \\
4x\phi(y-x^2)\phi'(y-x^2) & -2\phi(y-x^2)\phi'(y-x^2) & 0
\end{bmatrix}
$$
which represents $J(f)$ as $I+N$, with $N$ nilpotent
(it easy to verify, using a computer symbolic algebra program, that $N^3 = 0$). 
Note that $v-u^2-w = y-x^2$.
Put $a = \phi(y-x^2) = \phi(v-u^2-w)$. Then $f(x,y,z) = (x+a,y+z+2xa,z-a^2)$
and one easily obtains the inverse of $f$ by solving
for $z,x,$ and $y$, in that order. The result is
$$
(x,y,z) = f^{-1}(u,v,w) = (u-a,v-w-2ua+a^2,w+a^2)
$$
\end{Ex}

These examples are chosen to illustrate several points.

In Example \ref{Ex1}, it is clear that $3$ and $5$ could
be replaced by any other pair of constants, and that the cosine
could be any $C^1$ function, and the example would still
have the same properties. In fact, the example is perfectly
general, since any $C^1$ map of $\R^2$ to itself
with a unipotent Jacobian matrix has,
up to a constant translation, this form for some pair
of constants and some $C^1$ function. A proof is supplied
later.  

Example \ref{Ex2} shows how tacking on additional
coordinates in a triangular fashion allows them to parameterize
examples with fewer variables (if $d=0$, then $w$ is a parameter in the family
of $3$-variable triangular maps given by the first $3$ components).
In addition it demonstrates the existence of an explicit,
closed form inverse obtained by composition of functions.
The inverse is explicit in the sense that it is represented
by a (finite) formula involving algebraic operations, composition,
and the functions $a,b,c$ and the constant $d$ that appear in
the definition of $f$.

Example \ref{Ex3} shows how the principles involved in the first two
examples can be combined. It takes the first example, adds a
coordinate $z$ in triangular fashion, which therefore can
be used as a parameter in the first example, and then it
replaces $3$ by $1$, $5$ by $z$, and $\cos$ by $\phi$.
Engelbert Hubbers and Arno van den Essen introduced this sort of
bootstrap construction; 
if one adds as a construction technique the 
replacement of a map $f: \R^n \To \R^n$ by
the map $p \mapsto T^{a}f(Tp)$, where $T$ is a matrix of constants and $T^{a}$ is
its classical adjoint, one is led to their New Class of Automorphisms
\cite{NewClass}.
Actually, they consider only polynomials, but over
an arbitrary commutative coefficient ring with unit, which is both
a more restricted and more general case than that of the $C^1$ maps that are
the focus of this paper. 
All automorphisms in the New Class have unipotent Jacobian matrices
and the special property that for $f(p) = p + h(p)$, the perturbation
portion, $h: \R^n \To \R^n$, of $f$ has a constant $n$-fold composition power;
that is, $h^{\circ n} = h \circ h \circ \cdots \circ h$ ($n$ times) is a constant map.
This example, with $\phi(x+zy) = -(x+zy)^2$, produced the first three
dimensional counterexample to the Markus-Yamabe conjecture \cite{PolyCxMYC}.
In detail, the ordinary differential equation $dp/dt = -f(p)$ in $\R^3$ has
an orbit that escapes to infinity in forward time, 
even though the eigenvalues of
$J(-f)$, which are all $-1$, 
obviously have strictly negative real parts at
every point. 
One such orbit is $p(t) = (x(t),y(t),z(t)) = (18e^t,-12e^{2t},e^{-t})$.

Example \ref{Ex4}  is a very modest generalization of the
particular case $a = \phi(y-x^2) = y-x^2$ which
was constructed in \cite{IndependentRows}, where it
is shown that the perturbation part of $f$,
namely $h(x,y,z) = (y-x^2, z + 2x(y-x^2), -(y-x^2)^2)$, has nonconstant composition
powers $h^{\circ n}$ for every $n>0$. Thus $f$ is not generally an automorphism
in the New Class, since that particular $f$ is not.

\section{Conjectures}
The Jacobian Conjecture is a significant unsolved problem.
A polynomial map with a global polynomial inverse 
has a Jacobian matrix whose determinant is a nonzero constant.
This follows from the chain rule and the fact that the product
of two polynomials is a nonzero constant if, and only if, each
factor is. The Jacobian Conjecture is that the
converse is true. 
It is sometimes known as Keller's Jacobian Conjecture, because
its first appearance in the literature appears to be 
\cite{Keller}, in which Keller proves the complex birational case.
A modern formulation (\cite{Survey,ArnoBook}) is

\begin{Conj}[The Jacobian Conjecture]
Let $k$ be a field of characteristic zero, and $f: k^n \To k^n$
a polynomial map. Then $f$ has a polynomial inverse if, and only
if, the Jacobian matrix of $f$ has a nonzero constant determinant.
\end{Conj}

The complex case is known to be universal \cite{Survey}; that is,
if the conjecture is true for $k = \C$, then it is true in general.
The real case ($k=\R$) implies the complex case (consider $\C^n$ as
$\R^{2n}$ and compare the Jacobian determinants). For any $k$, the
conjecture can be reduced to the consideration of maps of the form
$f(x) = x + (Ax)^3$, where $A$ is a matrix of constants and the
cube $(Ax)^3$ is computed component-wise \cite{Effective}.
If the Jacobian determinant of such a
cubic-linear map is constant, then it has a unipotent Jacobian matrix. 
The reduction of the general case to the cubic-linear case involves
introducing extra variables in general; that is, the invertibility of
a polynomial map in $n$ variables with a nonzero constant Jacobian
determinant is equivalent to the invertibility of a cubic-linear
map in $m$ variables, where $m$ is usually (significantly) larger
than $n$.
So the cubic-linear variant implies the general case only if
one considers all $n$. The status of the conjecture at this moment
is that the cubic-linear and related cases have been affirmatively
resolved for certain low dimensions ($n \le 7$ -- 
see \cite{HubbersThesis,Seven}), 
that a few special cases work (e.g. maps with at most quadratic
terms;
see \cite{Survey} for more),
but that the general
case, though true for any $k$ and $n=1$, is not known for
any $n > 1$ for even a single field of characteristic zero.

Questions about the existence of an inverse for polynomial maps
hinge on injectivity; for polynomial maps over $\C$ or $\R$,
injectivity implies surjectivity \cite{IntoOnto,GeomAlgReele}.
The non-vanishing of the Jacobian determinant in these two cases
implies that the map is locally an analytic isomorphism. In the
complex case, an injective polynomial map with a nowhere vanishing
(and hence constant) Jacobian determinant is a birational map, and
it has a global polynomial inverse \cite{Keller,Rudin}. In the
real case, a nowhere vanishing Jacobian determinant need not be
a constant; if $f: \R \To \R$ is defined by $f(x) = x + x^3$, then
$\det J(f) = 1 + 3x^2$ is nowhere vanishing, and $f$ is injective, so it
is a global real analytic homeomorphism, but its inverse is not
polynomial. Pinchuk described \cite{StrongRealJC} a class of 
polynomial maps
$f: \R^2 \To \R^2$ with nowhere vanishing but non-constant Jacobian
determinant that are not injective, thus refuting what was known variously 
as the {\it Strong Real Jacobian Conjecture}, or just
the {\it Real Jacobian Conjecture}.
So a nonzero constant Jacobian
determinant is a necessary hypothesis in the real case
of the Jacobian Conjecture
(and that is what the term {\it Real Jacobian Conjecture}
is now usually taken to mean).

In the real case one is naturally drawn to the question of
what can be said of more general maps (real analytic, $C^1$, etc.).
Global univalence (injectivity) of maps $f: \R^n \To \R^n$ is a
large, well-studied topic \cite{GlobU,GlobalDiffs}. There are numerous
conditions that can be imposed to obtain univalence, ranging
from general topological conditions (local homeomorphism + properness)
to ones more closely connected to the Jacobian matrix (positive definiteness
conditions, special matrix types, Hadamard's integral criterion).
A nonzero constant Jacobian determinant, by itself, certainly does
not suffice to guarantee univalence; that is, the straightforward
generalization of the Jacobian Conjecture to $C^1$ maps is false.
As a simple example of that, one can take the analytic map
$f: \R^2 \To \R^2$ given by 
$f(x,y) = (\sqrt{2}e^{x/2}\cos(ye^{-x}),\; \sqrt{2}e^{x/2}\sin(ye^{-x}))$;
it has Jacobian determinant $1$, but is not injective
(e.g. the image of $\{x=0\}$ is a circle). The example is taken from
Brian Coomes' paper \cite{MapInjectivity}; it is also mentioned
in \cite{RnPolyFlows,RJCSM,MountainPass}.
Constancy of the Jacobian determinant is a global condition
on the pointwise spectrum (set of eigenvalues) of the Jacobian
matrix (see \cite{NearlySpectral,EhresmannFibrations} for 
some results in that general category).
In that regard, a recent conjecture  \cite{MountainPass}
is worth highlighting

\begin{Conj}[Chamberland]
If $f: \R^n \To \R^n$ is $C^1$ and the eigenvalues of $J(f)$ are globally
bounded away from $0$,
then $f$ is injective.
\end{Conj}

This hypothesis covers the case of constant eigenvalues, and
thus of a unipotent (all eigenvalues are $1$) Jacobian matrix.
This suggests formulating

\begin{Conj}[$C^1$ Unipotent Jacobian Univalence Conjecture]
If $f: \R^n \To \R^n$ is $C^1$ and the matrix $J(f)$ is unipotent, then
$f$ is injective.
\end{Conj}

This conjecture is at least known to be true for $n=2$
(see section \ref{PlanarCase}), and, by reduction to
the cubic-linear case, is strong enough to imply the
truth of the Jacobian Conjecture if it is established
for all $n$. All the examples of $C^1$ maps with unipotent
Jacobian matrices that were presented earlier are, in fact,
$C^1$ automorphisms; that is, they are not just injective,
but they also are surjective, and hence have a global
$C^1$ inverse. Thus no counterexample has been presented
yet to the somewhat stronger conjecture in which injectivity is
replaced by bijectivity.

In the context of conjectures related to the eigenvalues of
the Jacobian matrix, there are a few additional ones that
have arisen in connection with the study of the Markus-Yamabe
conjecture and its discrete (function iteration) analogue
\cite{%
DiscreteMYP,%
PolyCxMYC,%
StabilityInjectivityJC%
}.

\begin{Conj}[$C^1$ Stability Conjecture]
If $f: \R^n \To \R^n$ is $C^1$ and the eigenvalues of $J(f)$ 
have strictly negative real part at every point,
then $f$ is injective.
\end{Conj}

Note that the $3$-dimensional polynomial counterexample to the
Markus-Yamabe conjecture \cite{PolyCxMYC} presented above (Example
\ref{Ex3}) has an orbit that escapes to infinity, but it is still
injective, viewed as a map from $\R^3$ to $\R^3$.

\begin{Conj}[$C^1$ Fixed Point Conjecture]
If $f: \R^n \To \R^n$ is $C^1$ with $f(0) = 0$, 
and the eigenvalues of $J(f)$ have absolute
value less than $1$ at every point, 
then $0$ is the unique fixed point of $f$.
\end{Conj}

The Stability Conjecture for polynomial maps implies the Jacobian
Conjecture, and the Fixed Point Conjecture for polynomial maps
is equivalent to the Jacobian Conjecture 
\cite{DiscreteMYP,PolyCxMYC,Conjs+Pblms}.
The discrete Markus-Yamabe question (DMYQ) was raised in \cite{DiscreteMYP}:
do the hypotheses of the Fixed Point Conjecture
imply global convergence of iterates 
(for any $x_0$, the sequence $x_{k+1} = f(x_k)$ converges to $0$)?
A rational counterexample to DMYQ for $n=2$ is presented in the same paper.
Polynomial counterexamples to DMYQ been constructed \cite{ChaoticPolyAuts}
for $n \ge 4$ and for $n \ge 3$ in \cite{PolyCxMYC},
but they do have $0$ as a unique fixed point.

A conjecture equivalent to the $C^1$ Unipotent Jacobian Univalence Conjecture is
the following special case of the Fixed Point Conjecture.
See section \ref{Fix} for a proof of the equivalence of the
conjectures.

\begin{Conj}[$C^1$ Multiple Fixed Point Conjecture]
If $f: \R^n \To \R^n$ is $C^1$ and has two distinct fixed points,
then $J(f)$ has a nonzero eigenvalue at some point.
\end{Conj}

Another way to state this is: if $J(f)$ is nilpotent, then $f$ has
at most one fixed point.
Viktor Kulikov poses this conjecture for polynomial maps over $\C$
in \cite{JCandNilpotentMaps},
and shows that if it is
true (for all $n$), then the Jacobian Conjecture
follows. 
The Jacobian Conjecture also follows from the real case,
in view of the equivalence of this conjecture and the
Unipotent Jacobian Univalence Conjecture.
It should be mentioned that
\cite{DiscreteMYP} poses a different conjecture about
(complex polynomial) maps with nilpotent Jacobian matrices and refers
to it as the `Nilpotent Conjecture'; namely, that
the rows of the Jacobian matrix of such a map are linearly
dependent over $\C$. A counterexample is given in
\cite{IndependentRows}.

\section{Goals}
The main goal of this paper is to place the conjectures 
relating unipotence of the Jacobian matrix
and univalence of maps firmly in the $C^1$ domain.
To that end, proofs of the properties of maps corresponding
to the first three introductory examples are supplied in
the $C^1$ case.
The results are extensions
of what was known in the polynomial case.

All of the conjectures of the previous section are perhaps too ambitious. 
Arno van den
Essen, certainly an authority on the Jacobian Conjecture, has expressed,
in person and in print, the opinion that it may be true for $n=2$, but seems
unlikely to be true in general. And the other conjectures all
imply the Jacobian Conjecture, 
even if
they are established only for polynomial maps (in all dimensions).

Apart from the fact that these conjectures appear to be
natural and interesting in the larger domain of $C^1$
maps, there is a secondary goal in introducing them.
The $C^1$ hypothesis gives them more room to be wrong, and
it seems that effort
should be devoted to finding specifically non-polynomial
counterexamples, in the hopes that such will clarify the polynomial
situation.
Counterexamples
may shed
more light on what distinguishes the various cases
in the regularity hierarchy of
polynomial, rational, semialgebraic, analytic, smooth,
or just $C^1$ maps. Even a failure to find counterexamples
may be of some help.

The final portion of the paper deals with some tractable special
cases:
linearizable maps, maps with bounded images, 
and polynomial maps with no zeros at infinity.
While there are more specific individual results, what ties
these all together is the fact that, in each case,
if the Jacobian matrix of the map is nilpotent, then the
map has a unique fixed point.
The final section discusses the applicability of some of the ideas
explored in this paper to broader contexts than $C^1$ maps.
\section{The Planar Case}
\label{PlanarCase}
\begin{Thm} Let $f: \R^2 \To \R^2$ be $C^1$. Then $J(f)$
is unipotent if $f$ is of the form $$f(x,y) = (x + b\phi(ax+by)+c,
y - a\phi(ax+by)+d)$$ for some constants $a,b,c,d \in \R$ and some
function $\phi$ of a single variable. If that is the
case, then $f$ has an explicit global inverse. Conversely, if
$f$ is $C^1$ and $J(f)$ is unipotent, then 
$f$ is of the form shown (for a $\phi$ that is $C^1$).
\end{Thm}

Some remarks are in order before the proof. The analogous result
for polynomial maps over a field $k$ of characteristic zero 
is in \cite{Survey}. 
The precise form described above is not spelled out, but is
implicit in the proof of \cite[Theorem 6.2]{Survey};
\cite[Corollary 6.3]{Survey} explicitly states that if $f: k^n \To k^n$
is written as $f(x) = x - h(x)$ and $J(h)^2 = 0$, then for $n=2$ or $n=3$
it follows that $f$ is a (composition) product of elementary automorphisms,
and hence invertible.
A proof for polynomial maps with coefficients in a
$\Q$-algebra that is a unique factorization domain appears in \cite{NewClass};
it takes the case of a field $k$ as its point of departure.
The proof below is modeled on a proof by Yu Qing Chen for
holomorphic maps of $\C^2$ to itself \cite{HolomorphicUnipotent}.
Marc Chamberland has also proved
the same result for real analytic maps of $\R^2$ to itself
\cite{AnalyticUnipotent}, by entirely different methods.

The remainder of this section is devoted to the proof of the above
theorem and some relevant observations.

By adding a constant translation, which does not affect $J(f)$,
one can assume that $f(0)=0$, and then it suffices to consider the case 
$c = d = \phi(0) = 0$.
If $f$ is of the form mentioned, then it easy to establish that
$\phi$ is $C^1$ and 
$J(f)$ is unipotent, and to compute the inverse of $f$ (see
Example \ref{Ex1}). 
This leaves only the converse portion of the theorem to prove.

So suppose that $J(f)$ is unipotent, and that $f$ is $C^1$.
Denote points in $\R^2$ by $z = (x,y)$, and let
$h(z) = f(z) - z$. Then $J(h)$ is nilpotent; that is, $J(h)^2 = 0$.
If $h$ is constant, then $h=0$, and the desired representation
exists. So assume that $h$ is not identically $0$.
If one can show that $f$ is of the desired form
then at least one of $a$ and $b$ is not zero, or $h$ would be identically
zero. But then it is clear that $f$ is $C^1$ if, and only if, $\phi$ is $C^1$.

The main goal of the proof is to establish that
\begin{equation}
\label{eq:0}
h(z + (J(h)(z))\zeta)=h(z)
\text{ for all }
z,\zeta \in \R^2.
\end{equation}
For less notational clutter, write $A = J(h)(z)$. Then the goal
is to show that $h(z+A\zeta)$ is independent of $\zeta$.

\begin{remark}
Yu Qing Chen actually establishes this result for analytic maps
$h$ in any number of variables with $J(h)^2 = 0$.
His proof involves computing the power series expansion for $h(z+A\zeta)$
for fixed $z$
at $\zeta=0$
recursively, showing that all but the constant terms are zero.
It has no obvious extension to $C^m$ maps (even for $m = \infty$).
\end{remark}

If $A=0$ then it is clear that $h(z+A\zeta)$ does not depend on $\zeta$.
So let
$z \in \R^2$ be a point such that $A = J(h)(z)$ is nonzero.
Let $h(z) = (h_1(z),h_2(z))$.
The Jacobian matrix $J(h)$ is of constant rank 1 (since
it is nonzero but nilpotent) in a neighborhood of $z$, so (by
a classic theorem on Jacobian matrices of constant rank)
the functions $h_1$ and $h_2$ are dependent, in the precise sense that one of
them can be written as a $C^1$ function of the other
in a neighborhood of $z$
(see, for example, \cite[\S 98]{Kowalewski}).
Assume, without loss of generality, that $h_2(x,y) = g(h_1(x,y))$, where
$g$ is a $C^1$ function of one variable.
Since $J(h)^2=0$, its trace is zero, which yields 
$\partial h_1/\partial x + \partial h_2/\partial y =
 \partial h_1/\partial x + g'(h_1(x,y)) \partial h_1/\partial y = 0.$
But then the gradient of $h_1$ is $\partial h_1/\partial y$ times
the vector $(-g'(h_1(x,y)),1)$, which is constant along level curves of $h_1$.
This implies that the level curves of $h_1$ are straight line segments
(locally). Furthermore, the slope of the line segment which is the
level curve through a point $(x,y)$ is $g'(h_1(x,y))$, which is a
continuous function of $(x,y)$.

Now make an affine
change of coordinates, so that $z=0$ and the integral 
curve of $h_1$ through $z$
is horizontal. 
Then $\partial h_1/\partial x|_{(0,0)} =0$ and hence
$\partial h_1/\partial y|_{(0,0)} \ne 0$.
Consider (in the new coordinates) a small segment of the
$y$-axis, $|y| \le \eta > 0$, 
along which $\partial h_1/\partial y$ does not vanish.
Let the integral curve through $(0,y)$ have slope $\sigma(y)$.
Then $\sigma(0)=0$, and
$h_1(t,\sigma(y)t+y)=h_1(0,y)$ for $|t|$ small and $|y|$ small.
It follows that if $\eta$ is small enough, one can assume
that $\sigma(y)$ is well defined and $|\sigma(y)| \le 1$ for $|y| \le \eta$,
and that there
exists an $\epsilon > 0$, such that
\begin{equation}
\label{eq:1}
h_1(t,\sigma(y)t+y)=h_1(0,y) \text{ for all }
|t| \le \epsilon, |y| \le \eta.
\end{equation}

Suppose that $\partial h_1/\partial y$ has a value of $0$ at
a point $(t,\sigma(y)t+y)$ with $|t|=\epsilon$ and $|y| \le \eta$.
Then there is a sequence of difference quotients $\Delta h_1/\Delta y$
with limit $0$, computed using the fixed point $(t,\sigma(y)t+y)$ and a variable
point $(t,\sigma(y)t + y + \Delta y)$, and one can 
assume that the variable point is also the endpoint of a level curve of $h_1$
through a point $(0,y')$ with $|y'| \le \eta$ 
(this may involve considering one-sided differences only, if $|y| = \eta$).
Connect the fixed point and the variable point back to the
corresponding points $(0,y)$ and $(0,y') = (0,y+\Delta' y)$ on the $y$-axis
by moving along the (straight) level curves. Compute a difference
quotient $\Delta' h_1/\Delta' y$ for those points. By the invariance
of $h_1$ along the level curves $\Delta' h_1 = \Delta h_1$. 
Since $h_1$ is monotone on the segment
of the $y$-axis considered, 
$\Delta' y \ne 0$ and it has the same sign as $\Delta y$.
By considering the quadrilateral with
the 4 points in question as vertices, it is clear that the ratio
between $\Delta y$ and $\Delta' y$ is bounded above and
below by a function of the slopes of the level curves, and hence by
absolute constants for fixed $\epsilon$ and $\eta$.
This means that the difference
quotients $\Delta' h_1/\Delta' y$ tend to zero, which contradicts
the fact that $\partial h_1/\partial y $ does not vanish at $(0,y)$ for
$|y| \le \eta$.

The above argument establishes that the endpoints of the level
curves at $|t|=\epsilon$ are nonsingular points (ones at which
$J(h)$ does not vanish). The original argument then shows that
both $h_1$ and $h_2$ are constant on line segments near each 
endpoint, which implies that the identity $h_1(t,\sigma(y)t+y)=h_1(0,y)$
can be extended past each endpoint. By compactness, $\epsilon$ can
be made uniformly larger (keeping the same fixed $\eta$) in equation
\ref{eq:1}.
Since there is no maximum value of $\epsilon$ for which equation
\ref{eq:1} holds, all the integral curves in question must be entire
straight lines.

\begin{remark}
The argument to show that there is no maximum value for $\epsilon$
was suggested  (in a somewhat different form) by Yu Qing Chen 
(private correspondence).
\end{remark}

By considering $h_2$ instead of $h_1$, if necessary, it follows that
there is a straight line through any point $z$ at which $A=J(h)(z) \ne 0$,
such that both $h_1$ and $h_2$ are constant on that straight line.
Suppose that the straight line in question is $z+tw$, 
where $t$ is a parameter and $w$ is a fixed nonzero vector.
Differentiating with respect to $t$ yields $Aw=0$.
But $A$ is of rank $1$, so its nullspace equals its image. Thus
$A\zeta$ is a multiple of $w$ for any $\zeta \in \R^2$.
This establishes equation \ref{eq:0}. 
From here on, the argument is very similar to
that of \cite{HolomorphicUnipotent}.

Any two such lines through different points $z$ where $J(h)(z) \ne 0$
must be parallel.
For if not, they intersect, any further such line intersects at least
one of the two, and $h$ is constant on the union of all those lines.
But then $h$ is constant on (at least) the open set where $J(h)$
is nonzero; which is absurd. 
Now fix one $z$ with $A = J(h)(z) \ne 0$, and
consider the line as above through it.
Then $h$ is constant on every line parallel to that one (either because
$J(h)$ is $0$ at every point of such a line, or because it is not).

Next make an invertible linear change of coordinates (a rotation, for
instance) that makes the parallel lines vertical. In the new coordinate
system $J(h)$ still satisfies $J(h)^2 = 0$ (by similarity), and $h$
is a function of $x$ alone. Suppose $h(x) = (r(x),s(x))$. From $J(h)^2 = 0$,
one obtains $(r')^2 = 0$, hence $r'=0$, hence $r=0$ (since $r(0)=0$).
So $f$ is the triangular map $f(x,y) = (x,y+s(x))$. Changing coordinates
back to the original system yields the desired representation.

\begin{Ex}
The normal form given above exists, in general, only for maps
that are defined (and $C^1$) on all of $\R^2$.
The map $h(x,y) = (x/y, \ln(x/y))$ from $\{x>0,y>0\}$ to $\R^2$
has a nilpotent Jacobian matrix, but the level curves 
(of either component of $h$)
are rays
from the (excluded) origin into the first quadrant.
\end{Ex}
\section{Fixed Points, Inverses, and Composition Powers}
\label{Fix}
Let $f: \R^n \To \R^n$ be any map. For $a \in \R^n$, 
define $g_a(x)$ as $x-f(x)+a$.

\begin{Lem}
The family $g_a$ consists of maps with at most one
fixed point if, and only if, the map $f$ is injective. 
\end{Lem}

\begin{proof}
Let $p$ and $q$ be two distinct fixed points of $g_a$. 
Then $f(p) = a = f(q)$. 
And conversely.
\end{proof}

The unipotence of $J(f)$ is equivalent to the nilpotence of
any or all of the maps $g_a$.
So the $C^1$ Unipotent Jacobian Univalence Conjecture
and the $C^1$ Multiple Fixed Point Conjecture
are equivalent, and even separately so for each dimension $n$ and
particular map $f$ and family $g_a$.

\begin{remark}
Call a $C^1$ map $g$ non-degenerate at a point $p$ if $J(g)(p)$
does not have $1$ as an eigenvalue.
If $p$ is a non-degenerate 
fixed point of $g$, then $f(p) = 0$, where $f(x) = x - g(x)$.
Since $J(f)(p)$ has no zero eigenvalue, it is invertible,
and $f=0$ has a unique solution in a neighborhood of $p$.
This shows that non-degenerate fixed points are isolated.
In particular, any fixed points of a
map $g$ with nilpotent $J(g)$ are isolated.
So the Multiple Fixed Point Conjecture holds for maps
with a fixed point set that is not discrete, and also
for maps with a connected fixed point set.
\end{remark}

Given $y = f(x) = x - g(x)$, it is natural to try to solve
for $x$, obtaining $x = y + g(x) = y + g(y + g(x)) = y + g(y + g(y + g(x)))
= \ldots$. This procedure terminates, and yields a the inverse
of $y$ under the map $f$, provided the chain of expressions eventually
becomes independent of $x$. The easiest way to see this is
the obvious

\begin{Lem}
$f^{-1}(y)$ consists of a single point $x$ if, and only if, $x$ is
the unique fixed point of the map $w \mapsto y + g(w)$.
\end{Lem}

As a corollary, one obtains
\begin{Thm}
\label{PowerThm}
Let $f: \R^n \To \R^n$ be a map, and $h: \R^n \To \R^n$ its
perturbation part, defined by $h(x) = f(x) - x$. For any
given point $y \in \R^n$,
suppose that the map $\tau_y(w) = y - h(w)$ has a composition power
$\tau_y^{\circ m} = \tau_y \circ \cdots \circ \tau_y$ ($m > 0$ times)
that is constant. Then $\tau_y$ has a unique fixed point $x$,
and $x = f^{-1}(y)$ can be computed by starting with any initial point and
applying $\tau_y$ $m$ times.
\end{Thm}

\begin{proof}
Let $p$ be the single point that is the image of $\R^n$ under
the constant map $\tau_y^{\circ m}$. Then $p$ is a fixed point
of $\tau_y$. Furthermore,
if $p$ and $q$ are fixed points of $\tau_y$, then
$p = \tau_y^{\circ m}(p) = \tau_y^{\circ m}(q) = q$.
\end{proof}

These are the operations that provide what were called
``explicit inverses'' in Examples 1--3. They involve
composition, but apart from that only algebraic operations
and the components of the map (granted, in terms of its
decomposition $f(x) = x - g(x)$). The simplest case is the
planar case. Given the representation
$f(x,y) = (x + b\phi(ax+by)+c,y - a\phi(ax+by)+d)$
of the previous section for $f(z) = f(x,y)$, one
can verify the composition power condition for $m=2$
and for any point whose inverse one wants to find.
In fact, as pointed out in \cite{HolomorphicUnipotent},
the inverse of $f$ can be written as $f^{-1}(z)
= z + g(z+g(z))$, since $z + g(z + g(w))$ is independent of $w$.
Alternatively,
$f^{-1}(z) =
z + g(z + g(-z)) = z + g(z -z -f(-z)) = z + g(-f(-z))$,
so
$f^{-1}(z) =
z + (-f(-z)) -f(-f(-z)) = z - f(-z) - f(-f(-z))$.
This shows that there is nothing unique about the way
of expressing the inverse in terms of the original map
and composition products.
Note that these expressions are universal formulas,
in the sense that they are valid for every $f$ with
unipotent $J(f)$, in exactly the same form.
Both representations can legitimately be considered
as explicit (even ``closed form'') formulas.
And they rely on a property, namely constancy
of the composition power, that is a nonlinear analogue
of nilpotence. The situation is reminiscent of the inversion
of a unipotent matrix $I - N$, where one has the ``explicit''
formula $I + N + N^2 + \ldots + N^m$ for the inverse when
$N^{m+1}$ vanishes.  Gary Meisters first drew attention
to this property of composition powers in \cite{InvertingPolyMaps}.

If one adds
the {\it unique fixed point operator} $\mu$ to one's repertoire,
the inverse can be expressed as $f^{-1}(y) = \mu_x(\tau_y(x))$,
even when finite composition products of the $\tau_y$ do not
all become constant for a fixed number of composition
factors. But that is akin to adding allowing an
infinitary operation; in practice, fixed points of maps are often found
by taking the limit of a sequence $x_{i+1} = h(x_i)$ of iterates.

But not always. For Example \ref{Ex4}, the composition powers $h^{\circ i}$ 
of the perturbation
part of the map $f(p) = p + h(p)$ are not always constant. 
The map $h$ does have a unique fixed point.
However, this does not imply that composition powers of such a map 
necessarily have a unique fixed point.

\begin{Ex}
\label{Ex5}
let $h(x,y,z) = (\phi(y-x^2),z+2x\phi(y-x^2),-\phi(y-x^2)^2)$ be the
perturbation part of the automorphism of Example \ref{Ex4}.
If $\phi(0) = 0$ and $\phi(-1) = -1$, then
$(0,0,0)$ is (easily shown to be) the unique fixed point of $h$, but
$h^{\circ 3} = h \circ h \circ h$ has a nonzero fixed
point, which is therefore a periodic point of $h$ of period $3$.
Specifically, $\{(-1,1,-1),(0,-1,0),(-1,0,-1)\}$ is an orbit.
So $h^{\circ 3}$ has at least $4$ distinct fixed points.
\end{Ex}
\section{Strong Nilpotence}
Strong nilpotence is a notion that was introduced in \cite{SNdim5} for
polynomial
maps $F$ from $\R^n$ to $\R^n$. The Jacobian matrix matrix $J(F)$ was
called strongly
nilpotent if all products of $n$ matrices of the form $J(F)(a_i)$,
for possibly
distinct points $a_i \in \R^n$, are zero.  
In \cite{StronglyNilpotent} the Jacobian matrix
of a polynomial
map $F: k^n \To k^n$ was defined to be strongly nilpotent if the matrix
product
$J(F)(a_1)J(F)(a_2)\cdots J(F)(a_n)$  is zero, where the $a_i$ are
distinct sequences of
independent variables of length $n$
(that is, ``generic points").  If $k$ is infinite, this is
equivalent to the
obvious generalization of the first definition to $k$. The following
somewhat more
general definition can be used in both these cases.

Call a family $\cal{F}$ of linear endomorphisms of a vector
space $V$ {\em strongly nilpotent} if  there is a positive integer $r$
such that the
product (i.e. composition) of any $r$ factors in the family is zero.

\begin{remark} The use of families (indexed collections) of linear
transformations,
rather than sets, is sanctioned by the traditional terminology of
``commuting
families of endomorphisms.'' It serves, perhaps, to emphasize that the
factors
in a product in the above definition need not be distinct.
\end{remark}

A family $\cal{F}$ satisfying the above condition for a particular $r$ is
said to be strongly nilpotent of index $r$.
The smallest positive integer $r$ for which $\cal{F}$ satisfies the
above
condition is called the exact index of nilpotence of the family
$\cal{F}$.
The following elementary theorem can be derived from well known results
on the
simultaneous linear triangularizability of matrices, such as McCoy's
Theorem (\cite{McCoy,MatrixAnalysis}). However, a simple proof,
valid over any field $k$, is presented below.  It is inspired by the
proof
of linear triangularizability for commuting families of endomorphisms
(a theorem of Frobenius \cite{Vertauschbare}) as presented in
\cite{Humphreys}.

\begin{Thm} 

\label{SNT}
Let $r$ be the exact index of nilpotence of a strongly
nilpotent family of endomorphisms of a vector space $V$ of  finite
dimension $n$. Then $r \le n$ and there is a basis of $V$ in which
all the members of the family are represented by strictly upper
triangular
matrices.
\end{Thm}

\begin{proof}
By induction on $n$. Let $A$ be a product of $r-1$ factors from the
family
$\cal{F}$ which is non-zero (the result being trivial if $r=1$). Let $W$
be the kernel
of $A$. $W$ is non-zero, of dimension less than $n$, and
invariant under $\cal{F}$, since $AB=0$ for any $B$ in $ \cal{F}$.  
By induction, there is a non-zero $w \in W$ annihilated by $\cal{F}$.
The family $\cal{F}$ acts on $V/wV$.
Lift a strictly  upper triangularizing basis of $V/wV$ to
$v_2,\ldots,v_n \in V$ 
and let $v_1 = w$.  
Let $V_0 = \{0\}$ and $V_i = kv_1 + \cdots + kv_i$ for $0 < i \le n$.
Then $\cal{F}V_i \subseteq V_{i-1}$ for $0 < i \le n$.
Thus $(v_1,\ldots,v_n)$ is a basis of $V$ in which all members of $\cal{F}$
are
strictly upper triangular.  Since the product of any $n$ strictly upper
triangular matrices is zero, it follows that $r \le n$.
\end{proof}
\section{Linear Triangularizability}
Arno van den Essen and Engelbert Hubbers 
showed \cite{StronglyNilpotent} that polynomial maps over an
arbitrary field which
are perturbations of the identity by a map with
strongly nilpotent Jacobian matrix 
are linearly
(unit upper) triangularizable.   
Their proof explicitly uses the fact that polynomials are
involved. The theorem below characterizes linearly triangularizable
$C^1$ maps in similar terms, and also in terms of properties of
composition products involving the perturbation part.

The following terminology will be used. A map $f: \R^n \To \R^n$ is (unit
upper)
triangular if  $f = (f_1,\ldots,f_n)$ and $f_i  - x_i$ depends only on
those $x_j$ with $j>i$.
The qualification ``unit upper" will be implicitly assumed, unless otherwise
stated. The map $f$ is linearly triangularizable if there is a
linear 
mapping $S$ from $\R^n$ to $\R^n$ which is invertible, and such that 
$S^{-1} \circ f 
\circ S$ is triangular. 
A map $\tau: \R^n \To \R^n$
will be called a translation-dilation if it has the form $\tau(x) = a +
\lambda x$, where
$a \in \R^n$ and $\lambda \in \R$.
All linearly triangularizable maps have unipotent Jacobian matrices
and, by Theorem \ref{PowerThm} and condition (3) below,
they have explicitly computable inverses.

\begin{Thm} Let $f(x) = x + h(x): \R^n \To \R^n$  be a  $C^1$ map.
Then the following conditions are equivalent:
\begin{enumerate}
\item $f$ is (unit upper) linearly triangularizable;
\item The family of all matrices $\{J(h)(x)\}$ (for $x \in \R^n$) is
strongly nilpotent;
\item Any composition product of maps from $\R^n$ to itself, in which $n$
or more
of the factors equal $h$ and the remaining factors are
translation-dilations, is constant.
\end{enumerate}
\end{Thm}

\begin{proof}

(1) $\Rightarrow$ (2). Suppose $f = S^{-1} \circ t \circ S$, with
$t$ triangular and with  $S(x) = Ax$ for an invertible
matrix $A$.
Then $J(h)(x) = A^{-1}(J(t)(S(x)))A - I$ and the
matrices $J(t)(S(x))$ are unit upper triangular.
Thus the family $\{J(h)(x)\}$, for $x \in \R^n$, is similar to
a family of strictly upper triangular matrices, and hence
strongly nilpotent.

(2) $\Rightarrow$ (3). Let $p = p_1 \circ \cdots \circ p_s$ be a
composition product
of maps $p_i: \R^n \To \R^n$. By the chain rule, $J(p) = (J(p_1)\circ q)
J(q)$, where
$q = p_2 \circ \cdots \circ p_s$. By induction, $J(p)(x)$ is a product
of  $s$ matrices
that are Jacobian matrices of the factors $p_i$ evaluated at various
points.
Factors that are translation-dilations yield scalar factors
(multiplication by $\lambda$)
in the matrix product, and can be moved to the front. If there are $n$
or more factors
$p_i$ equal to $h$, the matrix product is zero by strong nilpotence. But
$J(p) = 0$
implies that $p$ is constant.

(3) $\Rightarrow$ (1). Let $a_1,\ldots,a_n$ be $n$ points in $\R^n$.
Let $\tau_i(x) = a_i + \lambda_i x$ for non-zero $\lambda_i \in k$.
Consider the composition product 
$p = h \circ \tau_1 \circ h \circ \tau_2 \circ \cdots \circ h \circ \tau_n$.
By assumption it is constant. Its Jacobian matrix at a point $J(p)(x)$
is the product of $2n$ matrices, $n$ of which are scalar multiples of
the identity ($\lambda_i I$) and $n$ of which have the form $J(h)(a_i +
\lambda_i(z))$
for some point $z$. Since $J(p)(x) = 0$ and the product of factors
corresponding
to translation-dilations is a non-zero multiple of the identity matrix,
the product
of the remaining factors is zero. Now vary the $\lambda_i$ and
take a limit in which all $\lambda_i$ tend
to zero to obtain $J(h)(a_1)\cdots J(h)(a_n) = 0$ -- that is, strong
nilpotence.
\end{proof}
\section{The New Class of Automorphisms}
In \cite{NewClass,DNA} Arno van den Essen and Engelbert Hubbers
introduced a class of invertible polynomial maps more general than
linearly triangularizable maps, but with some of the same properties.
Maps in this {\em New Class} have been used to provide counterexamples
to a number of conjectures involving spectral conditions on the Jacobian
matrix and stability and linearizability of maps
\cite{PolyCxMYC,CxMeistersConj,ChaoticPolyAuts}. 

The invertible polynomial ``maps'' are defined for any coefficient ring
$A$ that is a commutative ring with $1$.  Because the coefficient ring
can be, say, a finite field, the maps are actually polynomial morphisms
of affine $n$-space over $A$, 
many of which may correspond to the same underlying set-theoretic map
$A^n \To A^n$.
A polynomial morphism of affine $n$-space to itself can, 
of course, be identified 
with an $n$-tuple $f=(f_1,\ldots,f_n)$ of polynomials in $n$ variables with
coefficients in $A$.

They are perturbations of the identity map; that
is,
each is of the form $f(x) = x + h(x)$, where $f$ and $h$ are 
$n$-tuples of polynomial in $n$ variables.
Of course, that is no real restriction, since any polynomial map
can be written as $f(x) = x + h(x)$, by simpling defining $h(x)$
as $f(x) - x$.

They are singled out from all polynomial maps by specifying what perturbations $h$ are permitted,
and the definition of the allowed perturbations is a recursive
definition in terms of the coefficient ring $A$. That is, to
define the allowed perturbations $n$ variables, one starts with
the allowed perturbations in $n-1$ variables, but over the
coefficient ring $A[x_n]$, and specifies what operations can
be applied to these to yield allowed perturbations in $n$ variables.

The primary focus of this paper is on $C^1$ maps, and the definition
in \cite{NewClass} is not applicable, because the relationship
between functions in $n$ variables and
those in $n-1$ variables cannot be described by a simple construction
that adds
a variable (such as forming polynomials or power series in an
additional variable), but only by fixing a variable.  

An appropriate (still inductive) definition of
allowed perturbations
in $n$ variables is the following.
The induction specifies increasingly larger subsets
$\cal{H}_{n,i}$ of the set of all $C^1$ maps from
$\R^n \To \R^n$, all of which consist of allowed
perturbations. The final, largest subset, $\cal{H}_{n,n}$, is also
denoted $\cal{H}_n$, and a $C^1$ map $f: \R^n \To \R^n$
will be said to be in the New Class if (and only if)
it has the form $f(x) = x + h(x)$, where $h \in \cal{H}_n$.
The definition uses the notion of a {\it parameter} of a
map $f$, which refers to a coordinate variable $x_i$,
for which $f_i = x_i$.
Obviously, if $f(x) = x + h(x)$ and $x_i$ is a parameter
of $f$, then $h_i = 0$.

\begin{definition}
Let $n$ be a positive integer, let $x=(x_1,\ldots,x_n)$ be
a fixed (linear) coordinate system for $\R^n$. Let $h: \R^n \To \R^n$
be a $C^1$ map, with components given by $h=(h_1,\ldots,h_n)$. Then define
\begin{enumerate}
\item $h \in \cal{H}_{n,1}$ if $h_2 = h_3 = \cdots = h_n = 0$
and $h_1$ is independent of $x_1$;
\item $h \in \cal{H}_{n,i}$, for $i > 1$, if  
$h_{i+1} = h_{i+2} = \cdots = h_n = 0$ and there exist a map
$\widetilde{h} \in \cal{H}_{n,i-1}$, a map $T: \R^n \To \R^n$
that is $C^1$ and linear in the first $i$ coordinates with
$x_{i+1},\ldots,x_n$ as parameters, and a map $C: \R^n \To \R^n$
whose components depend only on the parameters 
$x_{i+1},\ldots,x_n$, such that
\begin{equation}
  h(x) = T^a \circ \widetilde{h} \circ T + C \tag{*} 
\end{equation}
\end{enumerate}
\end{definition}

In equation (*) above, $T^a$ refers to the classical adjoint of $T$.
More precisely, since $T$ has the last $n-i$ variables as
parameters, and it is linear in the first $i$ variables,
it can be represented as (pre-)multiplication by an $n \times n$
matrix $M$ which is block diagonal, of the special form
$$
M = 
\begin{bmatrix}
B & 0 \\
0 & I
\end{bmatrix}
$$
with $B$ an $i \times i$ matrix, {\it whose coefficients
will in general depend on the parameters $x_{i+1},\ldots,x_n$}. 
Then $T^a$ denotes the map of the same form given by
replacing $M$ by its classical adjoint (matrix of cofactors).
Note also that the last $n-i$ components of $C$ are necessarily
$0$ if equation (*) is to be satisfied.

The above definition accurately captures the spirit of the
definition of the New Class of automorphisms in the polynomial
case. And the New Class has the properties
specified in Theorem \ref{NewClassThm}, below, which are
analogous to ones established in the polynomial case.
The proofs supplied, however, are different, as \cite{NewClass}
relies on notions of specialization and localization
appropriate only to the polynomial case.
The actual proof of the theorem
(which is more of an elaborate verification by induction)
is carried out in Appendix A.

\begin{Thm}
\label{NewClassThm}
Let $f(x) = x + h(x)$ be a map in the New Class; that is,
$h \in \cal{H}_n$.  Then
\begin{enumerate}
\item If $r \in \R$, then $rh \in \cal{H}_n$.
\item If $T$ is a linear map of $\R^n$ to itself,
then $T^a \circ h \circ T \in \cal{H}_n$.
\item If $C$ is a constant vector of length $n$, then $h+C \in
\cal{H}_n$.
\item $J(h)$ is nilpotent 
(and thus $J(f)$ is unipotent).
\item The composition power $h^{\circ n} = h \circ h \circ \cdots \circ h$
($n$ factors) is constant.
\item The map $f$ is invertible, with an explicit inverse. 
\end{enumerate}
\end{Thm}

\begin{remark}
If $T$ is a change of basis map, then the expression of $f$ in
terms of  the new basis is $T \circ f  \circ T^{-1}$ = $x + T \circ h
\circ T^{-1}$.
From the first two
properties
above, $h \in \cal{H}_n$ also belongs to $\cal{H}_n$ defined in terms of the
new,
linearly related, coordinate system. So $\cal{H}_n$
is actually invariant under linear changes of coordinates.
\end{remark}

Somewhat more is known in the case of polynomials with coefficients
in a commutative ring $A$ with $1$. For one thing,
\cite{NewClass} proves that if $f$ is an automorphism in the
New Class then so is its inverse, whereas the above theorem
only establishes the existence of an ``explicit'' inverse.
Also, \cite{NewClass} defines a subtly different class of
allowable perturbations $\overline{\cal{H}_n}(A)$ that is, in general,
a strict superclass of $\cal{H}_n(A)$.
Unpublished work of Arno van den Essen \cite{barH} shows that
for a field $k$ of characteristic zero one has $\overline{\cal{H}}_n(k) =
\cal{H}_n(k)$ for $n < 5$, and that the equality is false for $n \ge 5$.
Finally, the following more refined properties of $\cal{H}_n(A)$ are proved in
\cite{DNA}.
The relevant definitions can be found there.

\begin{Thm}(\cite{DNA})
If $f(x) = x + h(x)$ and $h \in \cal{H}_n(A)$ then
\begin{enumerate}
\item $f$ is a finite composition product of automorphisms of the form
$\exp(D)$, for locally nilpotent derivations $D$ satisfying $D^2(x_i) =
0$, $1 \le i \le n$.
\item $f$ is a stably tame automorphism.
\end{enumerate}
\end{Thm}
\section{Linearizability}
It is instructive to consider the linear case of the Multiple
Fixed Point Conjecture. A linear map $L$ 
always has at least
one fixed point, namely $0$. 
If $L$ has a second, distinct
fixed point $x$, then the equation $x = Lx$ expresses the fact
that $x$ is an eigenvector of $L$ with the nonzero eigenvalue $1$.
Roughly the same argument can be applied to a linearizable map.

Recall that $f: \R^n \To \R^n$ is said to be linearizable if there
is an invertible map $\sigma: \R^n \To \R^n$ such that 
$\sigma^{-1} \circ f \circ \sigma = L$ is a linear map.
This is a special case of the notion of conjugate maps:
$f$ and $g$ are conjugate maps if there is an invertible
map $\sigma$, such that $\sigma^{-1} \circ f \circ \sigma = g$.
The definition of conjugacy implies that $f$ maps some set to itself,
and that $g$ also maps some (possibly different) set to itself.
Thus one can speak of the fixed points of $f$ and $g$, and it is
clear that they correspond to each other bijectively 
($f(\sigma(q)) = \sigma(q)$ if, and
only if, $g(q) = q$).
If $q$ and $p = \sigma(q)$ are corresponding fixed points
and $f,g$ and $\sigma$ are all $C^1$ maps between open subsets
of $\R^n$, then the chain rule yields $J(g)(q) = A^{-1}J(f)(p)A$,
where $A = J(\sigma)(q)$. So nilpotence at a fixed point is
preserved under this type of conjugacy.

A (globally) linearizable map has at least one
fixed point: the one corresponding to the fixed point $0$ of
the linear map. 
A local linearization is a conjugacy relation in which $f$ 
and a linear map $L$ have corresponding fixed points $p$ and $0$.
The existence of, possibly local, linearizations
of various degrees of regularity (topological, $C^1$, analytic, polynomial) 
around a fixed
point is a subject of considerable interest in a number of 
areas of mathematics and has a vast literature. It has captured
the interest of the polynomial automorphism community recently
as a possible new approach to the Jacobian Conjecture
\cite{Conjugation,GlobalAnalyticConjugation,GlobalAttractor,CxMeistersConj,CxCLLC,ArnoBook}.

The following theorem shows that linearizability of
maps with nilpotent Jacobian matrices would imply the
Multiple Fixed Point Conjecture (and hence the Jacobian
Conjecture). It can also be considered as a special
case in which the conjecture holds.

\begin{Thm}
Suppose that $f: \R^n \To \R^n$ is $C^1$ with nilpotent fixed points.
If $f$ is globally linearizable, then $f$ has
a unique fixed point.
\end{Thm}

\begin{proof}
Suppose that $f$ and $L$ are (globally) conjugate.
Because the fixed points of $f$ are nilpotent (that is,
$J(f)$ is nilpotent at each fixed point), the fixed
point set of $f$ consists of isolated points
(see remark in section \ref{Fix}), hence 
is countable.
If $L$ has $1$ as an eigenvalue, then its fixed point
set is uncountable (as it contains all multiples of an
appropriate eigenvector). Thus $L$ does not have $1$
as an eigenvalue, hence its fixed point set consists of $0$
alone, and thus $f$ has a unique fixed point.
\end{proof}

\begin{remark} Note that the linearizing map $\sigma$ is
not even required to be continuous for this result.
The hypothesis asserts the nilpotence of the
Jacobian matrix only at the fixed points of $f$, rather than
generally. But that is all that is needed for this application.
As far as the differentiability of $f$ is concerned, that also
need be assumed only in the neighborhood of the fixed points,
so that the inverse function theorem can be applied locally
to show that the fixed points are isolated.
\end{remark}
\section{Counting Fixed Points}
Let $f: X \To X$ be a continuous map of a topological space
to itself and suppose that the
rational homology of $X$ is finite. 
That is, the homology groups $H_i(X,\Q)$ with rational coefficients
are all finite dimensional rational vector spaces and
all but finitely many of them are zero. The map $f$ induces
an endomorphism $f_* : H_*(X,\Q) \To H_*(X,\Q)$; that is, for each
integer $i$ a linear map $f_{*,i}: H_i(X,\Q) \To H_i(X,\Q)$.
Define the Lefschetz number $L(f)$ of $f$ 
as the alternating sum of the traces of the induced
maps: $L(f) = \sum (-1)^{i} \text{Tr}(f_{*,i})$.
The definition appears to depend on the choice of homology
theory (e.g. singular homology) and on the coefficient field
$\Q$; as usual, it is in fact independent of those choices
for ``reasonable spaces'' (e.g. homeomorphic to a polyhedron)
and ``reasonable coefficients'' (a field containing the rationals).
By using homology with integral coefficients (mod torsion),
one can also show in that case that $L(f)$ is an integer.
Two homotopic maps $X \To X$ have the same Lefschetz number.

Now let $M$ be a compact, oriented $C^1$ manifold and $f: M \To M$
a $C^1$ map. 
A fixed point $x$ of $f$ is said to be non-degenerate if the tangent map
$df_x$ from $T_x(M)$ to itself does not have $1$ as an eigenvalue.
In terms of local coordinates near the fixed point, this amounts to
the statement that the Jacobian matrix of $f$ in those coordinates,
evaluated at $x$, does not have $1$ as an eigenvalue (equivalently,
it has $0$ as its unique fixed point).
The local Lefschetz number of $f$ 
at a non-degenerate fixed point $x$ is defined
to be $+1$ if the determinant of $I-df_x$ is positive, and $-1$
if it is negative.
If all the fixed points of $f$ are non-degenerate, then there are
only finitely many of them (because they are isolated), and $L(f)$
is the sum of the local Lefschetz numbers of $f$ at all of its fixed points.
Since the local Lefschetz numbers are all either $+1$ or $-1$,
it follows that $L(f)$ is an integer that algebraically
(that is, with due attention to sign) counts the fixed points of $f$.
The non-degeneracy of all fixed points is geometrically equivalent
to transversality of the intersection of the graph of $f$ and the
diagonal submanifold $\Delta = \{(m,m) | m \in M\} \subset M \times M$.
It can be given a purely topological definition, and the local
Lefschetz number can be defined for continuous maps that satisfy
the non-degeneracy condition, with the same result that $L(f)$ is
the sum of the local Lefschetz numbers.
All of these facts are standard topological fare
\cite{TopologyAndGeometry,ModernGeometry}.

Obviously, a fixed point for which $df_x$ is nilpotent will have
a local Lefschetz number of $+1$. So if $df_x$ is nilpotent at
each fixed point of $f$, then $L(f)$ is the actual number of
fixed points of $f$.

Unfortunately, $\R^n$ is not compact (though life would be less
interesting otherwise). To apply the above results, one has
to pass to a completion of $\R^n$ that is a compact orientable
manifold. One possibility is $S^n$, the $n$-dimensional sphere,
which is the $1$-point compactification of $\R^n$. 
That idea leads to the proof of the following

\begin{Thm} If $f : \R^n \To \R^n$ is $C^1$ with nilpotent
fixed points and the range of $f$ is bounded,
then $f$ has a unique fixed point.
\end{Thm}

\begin{proof} 
Since $f(\R^n)$ is bounded, so is its closure. Let $D$ be a disc
(closed ball) large enough so that the closure of $f(\R^n)$ is
contained in its interior. Let $D'$ be a concentric larger
disk. Let $A$ be the closed annulus between the $n-1$ spheres
$\partial D$ and $\partial D'$. Then $A$ is homeomorphic to
$S^{n-1} \times I$. The restriction of $f$ to $\partial D$ maps
the $n-1$ sphere $\partial D$ into the contractible space $D$
(because $D$ contains the entire image of $f$). So it is homotopic
to a constant map. Interpret the homotopy as a map of $A$ into $D$
which coincides with $f$ on $\partial D$ and is constant on $\partial D'$.
Let $p$ be the constant terminal value of the homotopy on $D'$.
Now define a map $F$ on $S^{n} = \R^n \cup \{\infty\}$ as follows:
it coincides with $f$ on $D$, it is the selected homotopy on $A$,
and it is the constant $p$ on $S^n - D'$.
By construction,
$F$ is continuous (and it could be made $C^1$ as well),
and $F(\infty) = p$. Since $F(S^n) \subset D$, the only fixed points
of $F$ are in $D$, hence they are those of $f$. 
The homology groups of $S^n$ with
integral coefficients are $H_0(S^n,\Z) = \Z$ and $H_n(S^n,\Z) = \Z$,
with all others $0$. 
The induced map $F_*$ on $H_0$ is the identity.
On $H_n$ it is multiplication
by the degree of $F$. 
Since $F$ factors through a map to $D$, which is contractible,
$F_{*,n}$ is $0$.
Therefore
the Lefschetz number of $F$ is $(-1)^0 * 1 + (-1)^n * 0 = 1$.
So $F$ has exactly one fixed point, and thus so does $f$.
\end{proof}

\begin{remark} 
Robert Brown notes (personal communication) that the proof can be shortened
by invoking the Lefschetz-Hopf Theorem, and provides the
reference \cite{Fixpunkt}.
\end{remark}
\section{Polynomial Maps with No Zeros at Infinity}
The result of the preceding section does not apply to
polynomial maps in any significant way, since the only
polynomial maps that have a bounded image are constant maps.
In contrast, this section describes results specifically
for (some) polynomial maps.

Let $f: \R^n \To \R^n$ be a polynomial map.
Let $d_i$ be the total degree (degree in all the variables jointly)
of the component $f_i$ of $f$. 
The algebraic degree of $f$ is defined as $d = \max_i d_i$. 
Each $f_i$ can be written as a sum of homogeneous
polynomials (forms). The form of highest degree ($= d_i$)
is called the leading form of $f_i$.
In \cite{RealCase} John
Randall showed that if the Jacobian determinant of $f$ vanishes
nowhere in $\R^n$ and the leading forms of the components of
$f$ have no common non-trivial zeros, then $f$ is a proper map
and hence a diffeomorphism of $\R^n$ onto $\R^n$.

Since each leading form is homogeneous and $d_i > 0$ (by the
condition on the Jacobian determinant), the leading forms
all vanish at $0$. The condition in Randall's Theorem is that
this is the only common zero of all the leading forms.

Normally, when one speaks of zeros at infinity of a polynomial
map, the reference is to common non-trivial zeros of the forms
of degree $d$ in each of the components. 
Then the condition that $f$ have no zeros at infinity is
equivalent to the statement that the rational map of
real projective $n$-space to itself given by
$$(x_0:x_1:\cdots:x_n) \mapsto (x_0^d: x_0^d f_1(x/x_0):
\cdots : x_0^d f_n(x/x_0))$$ (in homogeneous coordinates)
actually defines a global continuous map of real projective
$n$-space to itself that extends the polynomial mapping $f$.
If $d_i < d$ for some
$i$, then the form of degree $d$ in that component is identically
zero. So, if $f$ has no zeros at infinity in the usual sense,
it satisfies Randall's condition {\it a fortiori}.

\begin{Ex}
Let $f: \R^2 \rightarrow \R^2$ be defined by $f(x,y) = (x+y^3,y-x^3)$.
Then $\det J(f) = 1 + 9x^2y^2$ and $y^3$ and $x^3$ have only $(0,0)$ as
a common zero. Thus $f$ is a diffeomorphism.
\end{Ex}

\begin{remark}
If $f(x) = x + (Ax)^3$ is a map in cubic-linear form whose Jacobian
determinant is constant (equivalently, $\det J(f) = 1$), then $\det A = 0$,
and $f$ therefore does have a zero at infinity.
\end{remark}

The following application to maps with non-degenerate Jacobian matrix
is basically just Randall's result from a different point of view.

\begin{Thm} Let $f: \R^n \To \R^n$ be a polynomial map 
whose Jacobian matrix does not have $1$ as an eigenvalue
anywhere (e.g. $J(f)$ is nilpotent).
If the leading forms
of the components of $f$ have degrees greater than $1$ 
and no common non-trivial zero,
then $f$ has a unique fixed point.
\end{Thm}

\begin{proof}
Consider the map $g(x) = x - f(x)$. 
Then $g$ has no zero eigenvalues anywhere, and so the
Jacobian determinant of $g$ vanishes nowhere.
The leading forms of the components of $f$ and of $g$
are negatives of each other.
By Randall's Theorem, $g$ is a diffeomorphism, so
the equation $g(x) = 0$ has a unique solution,
which is also the unique solution to $x=f(x)$.
\end{proof}

For maps with no zeros at infinity, the statement
is not as complicated.

\begin{Thm}
Let $f: \R^n \To \R^n$ be a polynomial map
with no zeros at infinity. If $J(f)$ does not have $1$ as an eigenvalue
anywhere (e.g. $J(f)$ is nilpotent), then $f$ has a unique fixed point.
\end{Thm}

\begin{proof}
If the algebraic degree, $d$, of $f$ is zero, then $f$ is constant,
and hence it has a unique fixed point. If $d=1$, then the Jacobian
matrix of $f$ is constant, with no eigenvalue of $1$. Thus $g(x) = x - f(x)$
is an affine map with a nonsingular linear part, so it is bijective,
and hence $f$ has a unique fixed point. If $d > 1$,
the fact that $f$ has no zeros at infinity implies that the
same is true for $g$. Though $g$ may have some components $g_i$
of degree $d_i < d$, Randall's condition must be satisfied,
since the leading forms of the components of degree $d$ have no
common non-trivial zero. So $g$ is a diffeomorphism, and $f$ has
a unique fixed point.
\end{proof}
\section{A Note on Scope}
A uniform context of $C^1$ maps of $\R^n$ to itself was
adopted throughout for clarity. Clearly, differentiability
assumptions can be dropped in some results, where only the
group structure of $\R^n$ plays a role (e.g. Theorem
\ref{PowerThm}). At the other extreme, a lot of the results
can be extended to situations in which differentiability
is only assumed to exist in an abstract sense (as a derivation),
provided that it has adequate formal properties. Such properties
include being defined on all polynomial maps (and perhaps others),
being defined on compositions, and satisfying the chain rule.
However, closure (the derivative is differentiable) and an
inverse function theorem are not always needed. As a
substitute for continuity one can use an axiom to the effect
that constant is equivalent to finite-valued and to a vanishing
derivative. The sort of situations to which one can then apply
the results include, for example, holomorphic maps over $\C$,
rigid analytic maps over a $p$-adic field, 
$C^1$ semialgebraic maps over a real closed ground field,
and other classes of maps of  $k^n$ to
itself for $k$ a field of characteristic zero.
\appendix
\section{Proof of Theorem \ref{NewClassThm}}
The proof of Theorem \ref{NewClassThm}
consists of a sequence of lemmas and
propositions.
The following notation will be used.  For a given $1 \le i \le n$, let
$v$ denote the
``live'' variables and $z$ the parameters; $v=(x_1,\ldots,x_i)$ and $z$
is either the
empty sequence  (if $i=n$) or  $z=(x_{i+1},\ldots,x_n)$ (if $i < n$). All functions,
matrix entries,
and components will be $C^1$, and dependence on the variables and
parameters
will be denoted in the usual fashion (for example, $f(x_i,z)$ denotes a
function or map independent of
$x_1,\ldots,x_{i-1}$). Juxtaposition is used for matrix and
matrix-vector products. The
notation $M^a$ will be used for the classical adjoint of the matrix
$M$.  The inductive
definition of $h \in {\cal{H}}_{n,i}$ can then be stated as
follows: $h$ has the form
$S \circ \widetilde{h} \circ T$ with $\widetilde{h} \in
{\cal{H}}_{n,i-1}$, $T(v,z) = (M(z)v,z)$, and
$S(v,z) = (M^a(z)v+\eta(z),z)$ for an $i \times i$ matrix $M$ and an
$i$-vector $\eta$ of functions
that depend only on $z$ and not on $v$. In the inductive proofs that
follow the base case
($i=1$) is always obvious and will not be explicitly checked. 

\begin{Lem} Let $\delta(z),\eta(z)$ and $N(z)$ be, respectively, a
scalar, a vector of length $i$, and
an $i \times i$ matrix, consisting of functions
that depend only on
$z$. If $h$ belongs to
${\cal{H}}_{n,i}$, then so do $\delta(z)h, h + (\eta(z),0)$ and
$(N^a(z)v,z) \circ h \circ (N(z)v,z)$.
\end{Lem}
\begin{proof}
Obvious induction. The last claim uses $N^a M^a = (MN)^a$.
\end{proof}

This establishes the first three claims in the theorem when $i=n$.

\begin{Lem}
Suppose $h \in
{\cal{H}}_{n,i}$.
Then the $i$-th leading principal submatrix of $J(h)$ is nilpotent.
\end{Lem}

\begin{proof}
The matrix in question consists of the first $i$ rows and columns of
$J(h)$. Temporarily
fix the values of the parameters. With $z$ fixed, apply the chain rule
to see that the
matrix is $M^a N M$, where the leading principal submatrix of rank $i-1$
of $N$ is
nilpotent of index $i-1$ (by induction) and the last row of $N$ is
$0$. It follows that
$N^i = 0$ and hence the same is true for $(M^aNM)^i = M^a N M M^a N
\ldots M =
\det(M)^{i-1} M^a (N^i) M = 0$. Since that is true for every fixed value
of $z$, the
desired result holds.
\end{proof}

This establishes the fourth claim in the theorem when $i=n$.

\begin{Lem}
If $\delta(z)$ and $\epsilon(z)$ are a scalar function and an $i$-vector 
of functions that depend
only on $z$,
and $h \in {\cal{H}}_{n,i}$, then $h \circ
(\delta(z)v+\epsilon(z),z) \in {\cal{H}}_{n,i}$.
\end{Lem}

\begin{proof}
$h(v,z) \circ (\delta(z)v+\epsilon(z),z)  = (M^a(z)v+\eta(z),z) \circ
\widetilde{h}(v,z) \circ (M(z)v,z) \circ (\delta(z)v+\epsilon(z),z) =
(M^a(z)v+\eta(z),z) \circ \widetilde{h}(M(z)\delta(z)v + M(z)\epsilon(z),z)
= (M^a(z)v+\eta(z),z) \circ \breve{h}(v,z) \circ (M(z)v,z)$, where
$\breve{h}(v,z) = \widetilde{h}(\delta(z)v
+ \gamma(z),z)$ and $\gamma(z) = M(z)\epsilon(z)$. By induction, and the
fact that all the parameters
are also parameters for $i-1$, it follows $\breve{h} \in
{\cal{H}}_{n,i-1}$ and hence that $h$
has the desired form.
\end{proof}

\begin{Lem} If $h \in {\cal{H}}_{n,i}$ then the composition
power $h^{\circ i} = h \circ h \circ
\cdots \circ h$ ($i$ times) depends only on the parameters $z$.
\end{Lem}

\begin{proof}
Let $j=i-1$.
$h^{\circ i} = (S(v,z) \circ \widetilde{h}(v,z) \circ T(v,z))^{\circ i} = S(v,z) \circ
(\widetilde{h}(v,z) \circ T(v,z) \circ S(v,z))^{\circ j}
\circ \widetilde{h}(v,z) \circ T(v,z)$, where $T(v,z) = (M(z)v,z)$ and
$S(v,z) = (M^a(z)v+\eta(z),z)$. 
$T(v,z) \circ S(v,z) = (M(z)M^a(z)v + M(z)\eta(z),z) = (\delta(z)v +
\epsilon(z),z)$ with $\delta(z) = \det(M(z))$ and $\epsilon(z) =
M(z)\eta(z)$. By the previous lemma, and the fact that parameters for
$i$ are parameters for $j=i-1$, it follows that 
$\widetilde{h}(v,z) \circ T(v,z) \circ S(v,z) \in {\cal{H}}_{n,j}$.
By
induction, the result of raising this transformation to the 
composition power $j$
depends only on the
parameters for $j$ -- that is, on $x_i$ and $z$.  Denote the power by
$g(x_i,z)$. Then
$h^{\circ i} =
S(v,z) \circ g(x_i,z) \circ \widetilde{h}(v,z) \circ T(v,z)$. But $g(x_i,z)
\circ \widetilde{h}(v,z) = g(0,0)$
since all components of $\widetilde{h}$ are zero from the $i$-th on. So
$h^{\circ i}$ depends only on $z$.
\end{proof}

This establishes the fifth claim in the theorem when $i=n$.

\begin{Prop}
For $f$ as in the theorem, $f$ has an explicit
inverse.
Specifically, for each fixed $y$, the value of $f^{-1}(y)$ is the unique
value
of the composition power $(y-h(x))^{\circ n}$ (which is a constant map; that is,
independent of $x$). 
\end{Prop}

\begin{proof}
Let $y \in \R^n$. By  claims (1)--(3), the map $y-h(x)$ is in
${\cal{H}}_n$.
By claim (5), its $n$-th composition power is constant. Since this is
true for
every $y$, it follows from Theorem \ref{PowerThm} that $f^{-1}(y) =
(y-h(x))^{\circ n}$ (composition power).
\end{proof}

This establishes the sixth claim in the theorem.
\newcommand{\noopsort}[1]{}
\providecommand{\bysame}{\leavevmode\hbox to3em{\hrulefill}\thinspace}

\end{document}